\documentclass[10pt,conference]{IEEEtran}
\IEEEoverridecommandlockouts
\usepackage{cite}
\usepackage{amsmath,amssymb,amsfonts}
\usepackage{algorithmic}
\usepackage{graphicx}
\usepackage{textcomp}
\usepackage{xcolor}
\usepackage{ulem}
\usepackage{comment}
\usepackage{caption}
\usepackage{subcaption}
\def\BibTeX{{\rm B\kern-.05em{\sc i\kern-.025em b}\kern-.08em
    T\kern-.1667em\lower.7ex\hbox{E}\kern-.125emX}}

\colorlet{FigColor}{green!30!blue!70!black!100!}
\colorlet{darkBlue}{blue!30!black!70!}

\newcommand{\arr}[1]{\textcolor{red}{Razi: #1}}
\newcommand{\hs}[1]{\textcolor{olive}{#1}}

\def \myver{2}
\def \newver{2}

\begin{document}

\title{Diversity Maximized Scheduling in RoadSide Units for Traffic Monitoring Applications}

 \author{\IEEEauthorblockN{Ahmad Sarlak}
 \IEEEauthorblockA{\textit{School of Computing} \\
 \textit{Clemson University}\\
 Clemson, SC, USA \\
 asarlak@clemson.edu}
 \and
 \IEEEauthorblockN{Xiwen Chen}
 \IEEEauthorblockA{\textit{School of Computing} \\
 \textit{Clemson University}\\
 Clemson, SC, USA \\
 xiwenc@g.clemson.edu}
 \and 
 \IEEEauthorblockN{Rahul Amin}
 \IEEEauthorblockA{\textit{Tactical Networks Group} \\
 \textit{MIT Lincoln Laboratory}\\
  Lexington, MA, USA \\
 rahul.amin@ll.mit.edu}
 \and
 \IEEEauthorblockN{Abolfazl Razi}
 \IEEEauthorblockA{\textit{School of Computing} \\
 \textit{Clemson University}\\
 Clemson, SC, USA \\
 arazi@clemson.edu}}


\maketitle

\begin{abstract}
This paper develops an optimal data aggregation policy for learning-based traffic control systems based on imagery collected from Road Side Units (RSUs) under imperfect communications. Our focus is optimizing semantic information flow from RSUs to a nearby edge server or cloud-based processing units by maximizing data diversity based on the target machine learning application while taking into account heterogeneous channel conditions (e.g., delay, error rate) and constrained total transmission rate. As a proof-of-concept, we enforce fairness among class labels to increase data diversity for classification problems. The developed constrained optimization problem is non-convex. Hence it does not admit a closed-form solution, and the exhaustive search is NP-hard in the number of RSUs. To this end, we propose an approximate algorithm that applies a greedy interval-by-interval scheduling policy by selecting RSUs to transmit. We use coalition game formulation to maximize the overall added fairness by the selected RSUs in each transmission interval. 
Once, RSUs are selected, we employ a maximum uncertainty method to handpick data samples that contribute the most to the learning performance. 
Our method outperforms random selection, uniform selection, and pure network-based optimization methods (e.g., FedCS) in terms of the ultimate accuracy of the target learning application.

\end{abstract}

\section{Introduction}

RoadSide Units (RSUs) are an integral part of smart transportation systems due to their role in communicating with vehicles and collecting visual information to develop temporal and spatial traffic flow models. This information can be used to manage traffic flow, redesign traffic systems, analyze traffic safety, and guide drivers and police officers to react properly to temporal traffic events \cite{chen2022network}.
Particularly, with the recent advances in Vehicle-to-Vehicle/Infrastructure (V2V/V2I) communications in modern wireless systems (e.g., mode 3 and 4 side channel service in LTE release 14 \cite{toghi2018multiple},  \cite{sarlak2021approach}, New Radio NR-5G and WiFi-based IEEE 802.11p and newer IEEE 802.11bd \cite{naik2019ieee}), Autonomous Vehicles (AVs) can assist RSUs to collect massive traffic imagery data from multiple points of view, in addition to their own cameras.

With the rise of Edge Computing (EC), the bulk of heavy computations can be performed in edge servers located in RSUs or the entry point of the wireless network in the vicinity of RSUs. For instance, tasks like denoising, distortion removal, perspective transformation, video stabilization, and video-based action recognition and anomaly detection can be performed individually in RSU/EC servers \cite{lin2022low}. Even part of the AV computations can be offloaded to RSU/EC servers \cite{amjid2020vanet}. However, to develop learning-based systems, we often need to use data/model sharing between RSUs and cloud servers. For instance, traffic sign imagery can be collected from RSUs to develop a universal traffic sign recognition model in a centralized server \cite{xie2022efficient,mchergui2022survey,boroujeni2021data}. 
Another example is a recent work by Chen et. al. \cite{chen2022network} which offers a network-level traffic safety analysis system where video pre-processing is carried out in RSUs, and Deep Learning (DL) models are used in a central unit to extract network-level safety metrics for the overall safety profiling of the highway traffic. Some other applications include pedestrian detection \cite{ojala2019novel}, and unknown object detection \cite{chen2023roadside} based on data collected from RSUs. A complete review of such systems can be found in \cite{bai2022infrastructure,guerna2022roadside,razi2022deep}.

An important consideration in building distributed learning systems is investigating the capacity and limitations of the underlying wireless network that can influence the operational performance and quality of the learning-based analysis and decision-making platforms \cite{chen2020joint}. Particularly, errors and delays in an imperfect communication system, along with its limited bandwidth and constrained transmission resources \cite{chinipardaz2022inter}, can substantially compromise the quality of the trained models \cite{li2020federated,owfi2023meta}. Any reduction in throughput means collecting fewer samples which translates to a decline in the prediction power of Machine Learning (ML) algorithms for training-based models such as vehicle classification, traffic light recognition, etc. \cite{wang2019adaptive}\cite{tran2019federated}. Likewise, unexpected delays and connectivity issues can jeopardize the real-time operation of the system for applications like online tracking, anomaly detection, and accident risk analysis \cite{chen2020convergence}. High packet drop rates may decline the quality of distributed learning systems by disrupting model and data sharing \cite{chen2021distributed}. These limitations can cause catastrophic consequences for traffic control platforms that operate based on the imagery collected from RSUs \cite{wong2022virtual}. Cyber attacks and adversarial learning attacks can mislead the trained models and discredit their resulting decisions \cite{kaviani2022adversarial}.

A large body of work is devoted to characterizing the impacts of networking factors on the quality of Distributed Deep Learning (DDL) tasks as well as enhancing networking performance to improve the ultimate quality of DDL. 
For instance, in edge computing in the internet of vehicles based on learning applications, lowering delay by optimizing network resource utilization is addressed in \cite{zhang2022task, abad2020hierarchical}. 
In \cite{guo2022distributed}, the authors study the trade-off between delay and energy consumption by optimizing one factor when the other is constrained. They used a federated learning model to solve these optimization problems. 
In \cite{samarakoon2019distributed}, the problem of joint power and resource allocation for ultra-reliable low-latency communication in vehicular networks based on a federated learning approach is studied. 

Another closely-related line of research is adapting data-sharing and model-sharing strategies based on the limitations of the underlying network.  
Konecny et al. presented a novel algorithm for federated learning in which the communication cost is minimized \cite{konevcny2015federated}. Lotfi et al. proposed a novel semantic-aware Collaborative Deep Reinforcement Learning (CDRL) method that enables a group of heterogeneous untrained agents with semantically-linked DRL tasks to collaborate efficiently across a resource-constrained wireless cellular network \cite{lotfi2022semantic}.

We consider this issue from a substantially different perspective by regulating packet transmission under imperfect networking so that the \textit{diversity} of accumulated samples in the processing unit is maximized. More specifically, we aim to optimize data collection from a set of RSUs with heterogeneous networking conditions when the total data aggregation limit is constrained. Our goal is to optimize \textit{semantic information exchange}, which may not necessarily translate to optimizing the raw data throughput. This approach is driven by the fact that some data samples may not significantly contribute to the ultimate learning quality \cite{lin2017focal}. For instance, it is known that more balanced datasets with an almost equal number of data classes can typically yield better classification results compared to unbalance datasets of the same size \cite{kaur2019systematic}. A similar fact applies when data representation in some potentially hidden space is more diverse \cite{yu2020learning}. For instance, imagery collected from different RSUs may contain an unbalanced number of vehicle types for different traffic compositions. Therefore, building an object-tracking model based on data collected from RSUs that lack some specific vehicle types may perform poorly when tracking such vehicle types.

To this end, we optimize the semantic diversity of the collected data by enforcing \textit{fairness} among data classes. This determines the number of samples per class. To select samples within each class, we select samples that can contribute the most to the ML application at hand. One approach for multi-level classification is using min-margin (i.e. the difference between the softmax highest and second highest probability). 
Our approach is different than pure communication-based optimization methods, where the goal is optimizing communication performance metrics such as delay and throughput without considering the semantic content of data packets. 
Simultaneously, our perspective is different than typical fairness-imposing scheduling methods that try to balance resource utilization by different network nodes without considering their contribution to the performance of the target ML application \cite{chaieb2022deep}. We bridge these two perspectives by regulating packet transmission rate in a real-time fashion so that the performance metrics, along with the added fairness in each interval, remain maximum while obeying the transmission constraints. Our contribution can be summarized as i) optimized scheduling for imperfect communications while maintaining high diversity by imposing fairness among accumulated class labels, ii) using coalition game theory to characterize the added data diversity by any selection of RSUs in an interval-by-interval fashion, and iii) using constrained satisfaction problem to translate optimal attempt probabilities to a binary scheduling matrix.  

\section{System Model}

\begin{figure}[!ht]
\begin{center}
\setlength\belowcaptionskip{-2\baselineskip}
\centerline{\includegraphics[width=1\columnwidth]{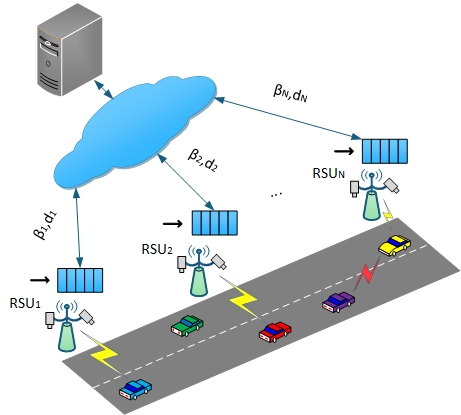}}
\caption{Roadside Units collect data samples from vehicles and send them to a central cloud-based computation server for learning-based processing.}
\label{fig:AccVsAPR}
\end{center}
\end{figure} 

We assume that there are $N$ RSUs represented by $n_1,n_2,\dots,n_N$ in a specific traffic zone. Each RSU, equipped with an Edge Computing (EC) server, pre-processes the collected imagery (e.g., video stabilization, projection transformation, etc.), then exploits training samples and sends them to a central processing unit through wireless links. There exist $M<N$ channels with an equal bandwidth $b_1=b_2=...b_M=bw$, so the number of simultaneous packets by all $N$ RSUs can not exceed $M$. We consider a slotted communication system with common channel access, a First Come First Serve (FCFS) queuing system with an infinite buffer, zero error tolerance protected by check codes like CRC, and selective auto-repeat request, so that all intended packets reach the destination error-free with random delays due to multiple attempts. Packet collisions are controlled by the proposed coordinated scheduling, so neither preventive techniques such as Carrier Sense Multiple Access (CSMA) nor collision management techniques such as Aloha are required.

\begin{table}[htbp]
\caption{List of Notations}
\begin{center}
\begin{tabular}{|c|l|}
\hline
\cline{1-2} 
\textbf{\textit{Notation}}& \textbf{\textit{Description}} \\
    \hline
    $N$ & Number of RSUs \\
    \hline
    $M$ & Number of available channels \\
    \hline
    $K$ & Number of active RSUs \\
    \hline
    $T$ & Number of timeslots in a transmission cycle \\
    \hline
    $n_i$ & the $i^{th}$ RSU \\
    \hline
    $\alpha_i$ & attempt probability of RSU $n_i$\\
    \hline
    $R$ & Number of re-transmissions \\
    \hline
    $D_i$\label{eq} & Delay of RSU $n_i$\\
    \hline
    $\beta_i$ & Packet error rate of RSU $n_i$ \\
    \hline
    $1/\lambda_i$ & Average packet delay of RSU $n_i$ \\
    \hline
    $\zeta_i$ & Throughput of RSU $n_i$\\
    \hline
    $C$ & Number of classes in the dataset \\
    \hline
    $c_i^j$ & Number of class $j$ in RSU $n_i$ \\
    \hline
\end{tabular}
\label{tab1}
\end{center}
\end{table}

To be more specific, suppose that a transmission cycle (we also call it interval) includes $T$ timeslots. We present the transmission matrix as a $Q_{N\times T}=[Q_{it}]$, where $Q_{it}=1$ means that RSU $n_i$ sends a packet at timeslot $t$ and is silent otherwise.

\begin{equation}  \label{eq:Q}
  Q_{N\times T} =
  \left[ {\begin{array}{cccc}
    1 & 0 & \cdots & 1\\
    0 & 1& \cdots & 0\\
    \vdots & \vdots & \ddots & \vdots\\
    1 & 1 & \cdots & 0\\
  \end{array} } \right]
\end{equation}
The goal of the scheduling is to fill in the binary scheduling matrix (interval by interval) so that the number of simultaneous transmissions in each timeslot does not exceed $M$, meaning that
\begin{align}   \label{eq:M}
&\sum_{i=1}^N Q_{it} \leq M,   &\forall t \in T.
\end{align}

We can solve this problem in two sequential steps. First, we define $\alpha_i=\frac{1}{T} \sum_{t=1}^{T} Q_{it}$ for RSU $n_i$, which represents the transmission attempt probability in a stationary case (during one transmission cycle: $T$ timeslots). Then, the optimization would reduce to finding the attempt probability vector $\boldsymbol{\alpha}= (\alpha_1,\alpha_2,\cdot\cdot\cdot,\alpha_N)$ so that 
\begin{align}
\sum_{i=1}^N\alpha_i&= \sum_{i=1}^N \big( \frac{1}{T} \sum_{t=1}^T Q_{it} \big) \\
\nonumber
&= \frac{1}{T} \sum_{t=1}^T \big(   \sum_{i=1}^N  Q_{it} \big) \leq   \frac{1}{T} \sum_{t=1}^T  M = M, ~~~~~~~ \forall i \in N. 
\end{align}

\begin{figure}[!ht]
\begin{center}
\setlength\belowcaptionskip{-2\baselineskip}
\centerline{\includegraphics[width=1\columnwidth]{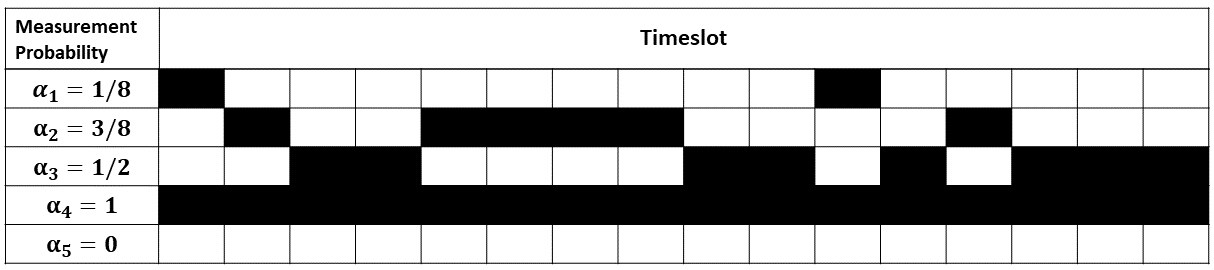}}
\caption{An exemplary pattern for N=5 targets and M=2 measurement resources.}
\label{fig:atp}
\end{center}
\end{figure}

Once, we obtain attempt probabilities $\alpha_i$ (the row weights of $Q$), then the next step would be to permute $\alpha_i T $ ones and $(1-\alpha)T$ zeros in row $i$ so that (\ref{eq:M}) is satisfied. We use constrained satisfaction problem by considering rows $r_i=(Q_{i1},Q_{i2},\dots,Q_{iT})$ as variables and all permutations of $\lfloor \alpha_i T \rfloor$ ones and $\lceil (1-\alpha_i)T \rceil$ zeros as values, where we used $\lfloor . \rfloor$ and $\lceil . \rceil$ to round up/down to an integer value, and select T to be the minimum product of $\lfloor \alpha_i M \rfloor$ terms. We will perform iterative assignment by random inconsistent variable selection and min-conflict heuristic for value selection until constraint (\ref{eq:M}) is satisfied. This will fully determine transmission matrix $Q$ with row weights $\alpha_i T$ and column vectors bounded by $M$. An exemplary scheduling for $N=8$, $M=2$, $T=16$, and $\mathbf{\alpha}=(1/8,3/8,1/2,1,0)$ is shown in Fig. \ref{fig:atp}.
A list of notations is given in Table I.

\section{Performance Metrics}
The problem is formulated as the following constrained optimization problem: 
\begin{align} \label{eq:opt}
\nonumber 
\arg\max _{\boldsymbol{\alpha}} f(\boldsymbol{\alpha},t)&= \sum_{i=1}^3 w_i f_i(\alpha_1,\alpha_2,\dots,\alpha_N,t) \\ \nonumber
&\text { s.t.}~~~~~ 0 \leq \alpha_i \leq 1 ,  \\  
&~~~~~~~~~~\sum_{i=1}^{\mathrm{N}} \alpha_i \leq \mathrm{M}, 
\end{align}
where the objective function is the linear combination of the set of desired performance metrics. In our case, $f_1(\boldsymbol{\alpha},t)$, $f_2(\boldsymbol{\alpha},t)$, and $f_3(\boldsymbol{\alpha},t)$, respectively, represent the average delay, throughput, and diversity metrics at time (interval) $t$. Here, the throughput and delay terms $f_1()$ and $f_2()$ are calculated per interval but the diversity term $f_3()$ depends on the RSU's transmitted packets in the previous timeslot, which hinders solving the problem in an interval-by-interval fashion. We will propose our game-theoretic strategy to address this matter in section IV.

The optimization is constrained by $0 \leq \alpha_i \leq 1$ meaning that each RSU can not utilize more than one channel at a time and $\sum_{i=1}^{\mathrm{N}} \alpha_i \leq \mathrm{M}$ meaning that the central unit can not receive more than $M$ packets simultaneously. To have a non-trivial solution, we set $M<N$. The components of objective functions are defined as follows.

\subsection{Delay} \label{sec:delay}

End-to-end (E2E) delay is an essential performance metric. Longer delays disrupt the online operation of time-sensitive tasks such as crash risk analysis. It also can compromise the accuracy of learning-based models and their adaptability to dynamic conditions. In delay-tolerant systems, longer delays can increase the packet timer expiration rate, leading to a higher packet drop rate. Although many efforts have been made to mitigate latency in modern communication systems, eliminating all delay sources is out of reach. E2E delays, in general, account for sampling and perception delays in the sender, queuing delays, channel setup delays, congestion and re-transmissions, actual transmission delays, as well as processing delays in the receiver \cite{cruz1991calculus}.    
We take a simplistic assumption and model $D_i^{(1)}$ the E2E delay for one packet for RSU $n_i$ as an exponentially distributed continuous-valued Random Variable (RV) with RSU-specific mean $\lambda_i$ following some prior work \cite{razi2017delay}, and then incorporate the impact of re-transmissions. Specifically, we define

\begin{align}   \label{eq:delay}
\nonumber
&D_i^{(1)}\sim f_{\lambda_i}(d_i), \\
&f_{\lambda_i}(d_i)= 
\begin{cases}
\lambda_i e^{-\lambda_i d_i}, & d_i\geq 0, \\ 0, & d_i<0.
\end{cases}
\end{align}

Here, we use capital letters for RVs and lowercase letters for their realizations. We consider that an average delay $E[D_i^{(1)}]=1 / \lambda_i$ for each RSU $n_i$ remains constant during one scheduling interval (one optimization round).   $\lambda_i$ can be drawn from Gamma distribution \cite{ali2008measured} in an interval-by-interval fashion.
Here, we consider $D_i^{(1)}$ captures all delay terms. However, if it accounts only for queuing delay, we can add the constant term $D_{tr}=L/ R_{ch}=L/ \alpha b_w$ accounting for the actual transmission delay, where $L$ is the packet length (in bits), $b_w$ is the bandwidth (in Hz) and $R_{ch}=\alpha b_w$ is the transmission rate (in bits/sec) for a given transmission system. 

This delay $D_i^{(1)}$ is for one packet transmission; hence, the effective E2E delay $D_i$ should be calculated for the last attempt when re-transmissions are involved. 
During one transmission interval, the number of re-transmissions $R_i$ for RSU $n_i$ follows a geometric distribution. Specifically, we have 
\begin{equation}
P (R_i=r)=(1-\beta_i)^{(r-1)}\cdot\beta_i,    \label{eq:geo}
\end{equation}
where $\beta_i$ is the packet drop probability for RSU $n_i$ drawn interval-by-interval from a Beta distribution with shape and scale parameters$\beta_i\sim \text{Beta}(a,b)$. Therefore, the expected value of E2E delay for one successful transmission is calculated as follows
\begin{align} \label{eq:dRSU1}
\nonumber
E\left[D_i\right]&=E_R \big [E_D[D_i \mid R\big] \\
\nonumber
&=E_R\big[E [ \underbrace{D_i^{(1)} +D_i^{(1)} +\dots D_i^{(1)}}_{R ~\text{times}} ] \big ]\\  
\nonumber
&=E_R\big[R \cdot   E [D_i^{(1)}] \big ]\\  
&=E[R] \cdot E[D_i^{(1)}]=\left(\frac{1}{\lambda_i}\right) \cdot\left(\frac{1}{1-\beta_i}\right),
\end{align}
where we used $E[R]=1.(1-\beta_i)$ and $E[D_i^{(1)}]=1/\lambda_i$.
The average delay of the entire system is apparently the linear combination of RSU-specific delays weighted by their attempt probabilities. We consider the negative/inverse of the average delay as our first term in the objective function in (\ref{eq:opt}). Therefore, we have \begin{align}  \label{eq:dtotal}
\nonumber
f_1(\boldsymbol{\alpha},t)&=1/ \sum_{i-1}^{N} \alpha_i(t) E [D_i (t)] \\
&=1/ \sum_{i-1}^{N} \frac{\alpha_i(t)}{(\lambda_i(t)) (1-\beta_i(t))} 
\end{align}
at transmission interval $t$. Note that we drop $t$ from some equations when it is clear from the context.

\subsection{Throughput} \label{sec:throughput}
Likewise, we can calculate the effective throughput $\zeta_i$ of RSU $n_i$ as  
\begin{equation}  \label{eq:th1}
    \zeta_i= \alpha_i R_{ch}/E[R_i] = \alpha_i R_{ch} (1 - \beta_i)
\end{equation}
where $R_{ch}$ is the rate of channel, $E[R_i]=1/(1 - \beta_i)$ is the average transmission per packet. 

The second term in the objective function of (\ref{eq:opt}) is the system throughput (at time $t$), which is simply the sum of the individual throughputs. We use
\begin{equation} \label{eq:th}
f_2(\boldsymbol{\alpha},t)= \sum_{i=1}^{\mathrm{N}} \zeta_i(t)=  \sum_{i=1}^{\mathrm{N}} \alpha_i(t)R_{ch} (1 - \beta_i(t)) .
\end{equation}

\subsection{Diversity through Fairness}
\ifx \myver \newver

Data collected from different RSUs may be extremely diverse due to different factors such as illumination conditions, camera resolution, camera altitude and field of view, background complexity, road geometry, and observed traffic composition. Therefore, enforcing fairness among RSU selection can enhance collected sample diversity.

It is known that the diversity of data samples significantly enhances learning quality. For instance, Determinantal Point Process (DPP) is used to increase the diversity of samples based on a kernel-based distance matrix \cite{kulesza2012determinantal}. DPP can be applied to the sample representation in the feature space or potentially lower dimensional representation space. Diversity can also be performed by graph analysis of data for more complex inter-sample dependencies \cite{mahmudi2019some,mahmudi2023some}. 
It is noteworthy that diversity may conflict with the efficiency of data accumulation (as shown in our results in section V). Specifically, an optimization merely based on delay and throughput would favor RSUs with better channel conditions, while diversity would favor more balanced scheduling. To address this trade-off, we enforce diversity through fairness among collected data samples of different classes. 
We can view it as \textit{soft fairness}.
There exist different fairness metrics, including min-max, alpha fairness, Jain’s index, and entropy. 
In this work, we use Jain's fairness index, defined as
\begin{align} \label{eq:fair}
J(\boldsymbol{x})=\frac{\big(\sum_{i=1}^{N} x_i\big)^2}{N \sum_{i=1}^{N} x_i^2 }.
\end{align}

To enhance diversity, we impose fairness among data categories (e.g., class labels in multi-level classification). Suppose $c_i^j$ is the number of data samples of class $j$ in RSU $n_i$. Then, the diversity of RSU $n_i$ is represented by the following vector
\begin{align}
\boldsymbol{c}_i=\left[c_i^1, c_i^2, \cdots, c_i^C\right]. 
\end{align}

Assuming that samples are selected at random by RSUs, then the number of received samples of class $j$ until time $t$, $c_R^j (t)$ is the sum of samples of the same class attempted by all RSUs in all transmission intervals up to time $t$ proportional to their effective throughput as follows
\begin{align}
c_R^j (t)= \sum_{\tau=0}^{\mathrm{t}}  \sum_{i=1}^{\mathrm{N}} c_i^j \zeta_i (\tau), \label{eq:fiar2}
\end{align}
where $\zeta_i(\tau)$ is the throughput of RSU $n_i$ at time $\tau$ given by (\ref{eq:th1}). It is clear that in contrast to delay $f_1()$ and throughput $f_2()$, the fairness component $f_3()$ can not be evaluated independently, because it depends on previous transmissions.
We impose fairness on the received data samples of all classes as

\begin{align}
{f_3(\boldsymbol{\alpha}, t)=} J\left(c_R^1(t), c_R^2(t), \cdots, c_R^C (t)\right). \label{eq:fiar3}
\end{align}

This completes the terms of the objective function in (\ref{eq:opt}). 

\else
Data collected from different RSUs may be extremely diverse due to different factors such as illumination conditions, camera resolution, camera altitude and field of view, background complexity, road geometry, and observed traffic composition.    
It is known that the diversity of data samples significantly enhances learning quality. For instance, Determinantal Point Process (DPP) is used to increase the diversity of samples based on a kernel-based distance matrix \cite{kulesza2012determinantal}. DPP can be applied to the sample representation in the feature space or potentially lower dimensional representation space. For instance, one may cluster data samples in lower dimensional space (e.g., based on an intermediate hidden layer with fewer neurons at the later stages of a convolutional neural network) and collects a balanced number of samples from the clusters. 

It is noteworthy that diversity may conflict with the efficiency of data accumulation (as shown in our results in section V). Specifically, an optimization merely based on delay and throughput would favor RSUs with better channel conditions, while diversity would favor more balanced scheduling. To address this trade-off, we enforce diversity through fairness among collected data samples of different classes. Fairness is often used to assign resources among users, but here we use it for a different objective of diversifying received data sample classes. We can view it as \textit{soft fairness}.

\hs{we can remove this part and just say that we use Jain's fairness index}
There exist different fairness metrics, including min-max, alpha fairness, Jain’s index, and entropy. An axiomatic view of fairness is given in \cite{lan2010axiomatic}, which combines a family of fairness metrics into the following unified metric:
\begin{align} \label{eq:fair}
f_\phi(\boldsymbol{x})=\operatorname{sign}(1-\phi) \cdot\left[\sum_{i=1}^N\left(\frac{x_i}{\sum_j x_j}\right)^{1-\phi}\right]^{\frac{1}{\phi}}
\end{align}
where sweeping $\phi$ from $-\infty$ to $\infty$ results in different fairness measures. For example, setting $\phi=-1$ and $\phi=0$, respectively, yields Jain's fairness index and Entropy. In this work, we arbitrarily use Jain's fairness index, defined as
\begin{align} \label{eq:fair}
J(\boldsymbol{x})=\frac{\big(\sum_{i=1}^{N} x_i\big)^2}{N \sum_{i=1}^{N} x_i^2 }.
\end{align}

To enhance diversity, we impose fairness among data categories (e.g., class labels in multi-level classification). Suppose $c_i^j$ is the number of data samples of class $j$ in RSU $n_i$. Then, the diversity of RSU $n_i$ is represented by the following vector
\begin{align}
\boldsymbol{c}_i=\left[c_i^1, c_i^2, \cdots, c_i^C\right]. 
\end{align}

Assuming that samples are selected at random by RSUs, then the number of received samples of class $j$ until time $t$, $c_R^j (t)$ is the sum of samples of the same class attempted by all RSUs in all transmission intervals up to time $t$ proportional to their effective throughput as follows
\begin{align}
c_R^j (t)= \sum_{\tau=0}^{\mathrm{t}}  \sum_{i=1}^{\mathrm{N}} c_i^j \zeta_i (\tau), \label{eq:fiar2}
\end{align}
where $\zeta_i(\tau)$ is the throughput of RSU $n_i$ at time $\tau$ given by (\ref{eq:th1}). It is clear that in contrast to delay $f_1()$ and throughput $f_2()$, the fairness component $f_3()$ can not be evaluated independently, because it depends on previous transmissions.
We impose fairness on the received data samples of all classes as

\begin{align}
{f_3(\boldsymbol{\alpha}, t)=} J\left(c_R^1(t), c_R^2(t), \cdots, c_R^C (t)\right). \label{eq:fiar3}
\end{align}

This completes the terms of the objective function in (\ref{eq:opt}). 

\fi

\section{Coalition-based Greedy Scheduling}
\ifx \myver \newver

The optimization problem in (\ref{eq:opt}) is non-convex, hence does not admit a closed-form solution or KKT approach. Solving numerically with an exhaustive search is computationally expensive. Specifically, if we divide the range of alpha $[0 ~~1]$ into $[0, d\alpha, 2d\alpha, \dots, 1]$ with $N_\alpha$ steps $d\alpha=1/N_\alpha$, the computation would be in the orders of $(N_\alpha)^N$ considering only O(1) complexity to evaluate objective functions $f_i(\boldsymbol{\alpha})$. Therefore, it is NP-hard in the number of RSUs and can be prohibitively expensive for large-scale systems.
The second and more important issue is that it does not allow interval-by-interval optimization, because the fairness index $f_3()$ should account for all transmission intervals to be evaluated at the end of transmission intervals. In other words, since the quality of channels is considered constant during one transmission interval, it is reasonable to evaluate performance metrics per interval but the fairness among accumulated class labels depends on the scheduling of preceding intervals.

To address this issue, we propose an approximate method to select top-$K$ RSUs using \textit{coalition game theory}. 
Game theory is an appropriate tool to evaluate fairness since it quantifies the contribution of each player $n$ when joining a coalition $S$, as \textit{marginal value} $v(S \cup \{ n \} )-v(S)$, where $v(S)$ is the value function representing the total payoff can be gained by the members of coalition $S$. Here we define 
\begin{align} \label{eq:opt2}
\nonumber 
v(S(t))= &\sum_{i=1}^3 w_i f_i(\alpha_1,\alpha_2,\dots,\alpha_N,t) \\ 
\nonumber 
 &w_1 / \sum_{n \in S(t)} \frac{\alpha_n (t)}{(\lambda_n(t)) (1-\beta_n(t))} \\
\nonumber 
 &+ w_2  \sum_{n \in S(t)} \alpha_n(t) R_{ch} (1 - \beta_n (t)) \\
&+ w_3 J\left(c_R^1(t), c_R^2(t), \cdots, c_R^C (t)\right)
\end{align}
with 
\begin{align}  \label{eq:equal-alpha}
\alpha_n = \begin{cases} M/|S(t)| & n \in S(t)   \\ 0 & else \end{cases}
\end{align}
which represents the case that members of coalition $S(t)$  (formed at transmission interval $t$) split $M$ transmission resources equally while the rest of $N-|S(t)|$ RSUs remain silent. 

In standard coalition games, we can have coalitions of arbitrary size. Also, Shapley value of each player is defined as the expected marginal contribution of player $i$ to the set of players who precede this player as
\begin{align}
\phi_i(v)= \frac{1}{N !}\sum_{S \subseteq \mathcal{N}/i}\frac{|S| !(N-|S|-1) !}{N !}[v(S \cup\{n_i\})-v(S)]
\end{align}
This involves evaluating the value of all $2^N$ coalitions. Here, to reduce complexity, we allow only the formation of fixed-size coalitions $|S(t)|=K$ as the set of active RSUs in transmission interval $t$. 
Apparently, we must have $|S|\geq M$ to maintain constraint $\alpha_i=M/|S|<1$. Note that $|S|$ can be much smaller than the number of RSUs, $|S|\ll N$. For instance, if we have $N=50$ and set $|S|=K=10$, then the number of coalitions reduces from $2^N\approx1.1\times 10^{15}$ to ${N \choose K}= \frac{50!}{10!(50-10)!}\approx 10^{10}$ (about 100,000 fold reduction).
This determines the set of active RSUs in each transmission interval.

\section{Uncertainty-Based Sample Selection }
Once we determine active RSUs by solving (\ref{eq:opt}) and (\ref{eq:opt2}), we can select a balanced number of samples among different classes. However, there is flexibility in selecting samples within each class. Our approach to this problem is sensing samples that can contribute the most to the ML application at hand. For instance, for multi-level classification we use the \textit{min-margin} criteria by selecting samples with the lowest difference between the softmax highest and second highest probability \cite{scheffer2001active}. 
To this end, the Fusion Center updates the model at the end of each transmission interval and sends back the model parameters. The RSU selects the samples that exhibit maximum uncertainty.

\else

The optimization problem in (\ref{eq:opt}) is non-convex, hence does not admit a closed-form solution or KKT approach. Solving numerically with an exhaustive search is computationally expensive. Specifically, if we divide the range of alpha $[0 ~~1]$ into $[0, d\alpha, 2d\alpha, \dots, 1]$ with $N_\alpha$ steps $d\alpha=1/N_\alpha$, the computation would be at least in the orders of $N_\alpha^N$ considering only O(1) complexity to evaluate objective functions $f_i(\boldsymbol{\alpha})$. Therefore, it is NP-hard in the number of RSUs and can be prohibitively expensive for large-scale systems.
The second and more important issue is that it does not allow interval-by-interval optimization, because the fairness index $f_3()$ should account for all transmission intervals to be evaluated at the end of transmission intervals. In other words, since the quality of channels is considered constant during one transmission interval, it is reasonable to evaluate performance metrics per interval but the fairness among accumulated class labels depends on the scheduling of preceding intervals.

To address this issue, we first select top-$K$ RSUs using \textit{coalition game theory}. 
Game theory is an appropriate tool to evaluate fairness since it quantifies the contribution of each player $n$ when joining a coalition $S$, as \textit{marginal value} $v(S \cup \{ n \} )-v(S)$, where $v(S)$ is the value function representing the total payoff can be gained by the members of coalition $S$. Here we define 
\begin{align} \label{eq:opt2}
\nonumber 
v(S(t))= &\sum_{i=1}^3 w_i f_i(\alpha_1,\alpha_2,\dots,\alpha_N,t) \\ 
\nonumber 
 &w_1 / \sum_{n \in S(t)} \frac{\alpha_n (t)}{(\lambda_n(t)) (1-\beta_n(t))} \\
\nonumber 
 &+ w_2  \sum_{n \in S(t)} \alpha_n(t) R_{ch} (1 - \beta_n (t)) \\
&+ w_3 J\left(c_R^1(t), c_R^2(t), \cdots, c_R^C (t)\right)
\end{align}
with 
\begin{align}
\alpha_n = \begin{cases} M/|S(t)| & n \in S(t)   \\ 0 & else \end{cases}
\end{align}
which represents the case that members of coalition $S(t)$  (formed at transmission interval $t$) split $M$ transmission resources equally while the rest of $N-|S(t)|$ RSUs remain silent. 

In standard coalition games, we can have coalitions of arbitrary size. Also, Shapley value of each player is defined as the expected marginal contribution of player $i$ to the set of players who precede this player as
\begin{align}
\phi_i(v)= \frac{1}{N !}\sum_{S \subseteq \mathcal{N}/i}\frac{|S| !(N-|S|-1) !}{N !}[v(S \cup\{n_i\})-v(S)]
\end{align}
which can be used to rank the users. \hs{we can remove this paragraph from the pape}. \hs{from here to the end }However, it involves evaluating the value of all $2^N$ coalitions. Here, to reduce complexity, we allow only the formation of fixed-size coalitions $|S(t)|=K$ as the set of active RSUs in transmission interval $t$. 
Apparently, we must have $|S|\geq M$ to maintain constraint $\alpha_i=M/|S|<1$. However, $|S|$ can be much smaller than the number of RSUs, $|S|\ll N$. For instance, if we have $N=50$ and set $|S|=K=10$, then the number of coalitions reduces from $2^N\approx1.1\times 10^{15}$ to ${N \choose K}= \frac{50!}{10!(50-10)!}\approx 10^{10}$ (about 100,000 fold reduction).

\arr{This is inconsistent with our assumption! Think how we can keep it more general with no restrictive assumption of $\alpha_i=...M/N$ for all active RSUs.}
Once, we identify the top-K RSUs, we can settle with equal transmissions among top-K RSUs or use numerical solutions to solve the optimization problem in (\ref{eq:opt}) for only top-K RSUs in a reasonable time. Our results show that these solutions are almost equivalent in long-run (for more than 100 transmission intervals).  \hs{end}

\section{Uncertainty-Based Sample Selection }
Once we identified the number of packets (or samples) for each RSU by solving (16) and developed the scheduling matrix $Q_{T\times N}$ in (1), the RSU must select the samples for transmission. Our approach to this problem is sensing samples that can contribute the most to the ML application at hand. One approach for multi-level classification is using 
\textit{margin} (i.e. the difference between the softmax probability of the highest and second highest predicted classes) \cite{scheffer2001active}. 
To this end, the [Fusion Center/ Central Unit] updates the model at then of each transmission interval and sends back the model parameters. The RSU selects the samples that exhibit maximum uncertainty.

\section{Uncertainty-Based Sample Selection }
Once we determine active RSUs, we can select a balanced number of samples among different classes.

identified the number of packets (or samples) for each RSU by solving (16) and developed the scheduling matrix $Q_{T\times N}$ in (1), the RSU must select the samples for transmission. Our approach to this problem is sensing samples that can contribute the most to the ML application at hand. One approach for multi-level classification is using 
\textit{margin} (i.e. the difference between the softmax probability of the highest and second highest predicted classes) \cite{scheffer2001active}. 
To this end, the [Fusion Center/ Central Unit] updates the model at then of each transmission interval and sends back the model parameters. The RSU selects the samples that exhibit maximum uncertainty.

\fi

\section{Simulation}

In this section, we investigate the performance of the proposed interval-by-interval scheduling policy in terms of the ultimate learning quality under time-varying conditions and compare it against uniform scheduling, random scheduling, as well as an exemplary communication-based scheduling method. In our experiments, we set the number of RSUs $N=10$, the number of available resources as $M=5$, the number of timeslots 
per each transmission interval to $T=10$, and the total number of intervals to $10$. Channel conditions ($\beta_i, \lambda_i$) are selected randomly and remain constant during one transmission interval. We execute the algorithm for $K=5$. We use Beta and exponential distributions to generate $\beta_i$ and $1/\lambda_i$, respectively.

We use the following two datasets for our experiments. The first dataset is \cite{jensen2016vision}, which provides $44000$ images of traffic scenes with traffic lights under different conditions taken by 5 RSU cameras, as shown in (Fig.2). We apply traffic light status detection (red, yellow, green) with a convolutional neural network (CNN) on this dataset as our exemplary application. Note that we only use a small subset of this dataset (about 150 images) for our test. 
We also evaluate our method by training a CNN on the benchmark CIFAR10 dataset \cite{krizhevsky2009learning} for multi-level classification. Our CNN model has 8 layers, the activation function is Relu, and the resolution of the input layer is $32 \times  32 \times 3$. Also, we use F1 score to measure the learning accuracy.

\begin{figure}[!ht] \label{fig:light}
\begin{center}
\setlength\belowcaptionskip{-2\baselineskip}
\centerline{
\includegraphics[width=1\columnwidth]{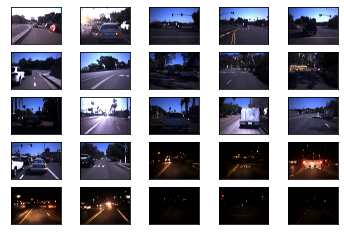}} 
\caption{Traffic light dataset. }
\label{fig:th}
\end{center}
\end{figure}

\begin{figure}[!ht]
\begin{center}
\setlength\belowcaptionskip{-2\baselineskip}
\centerline{
\includegraphics[width=0.95\columnwidth]{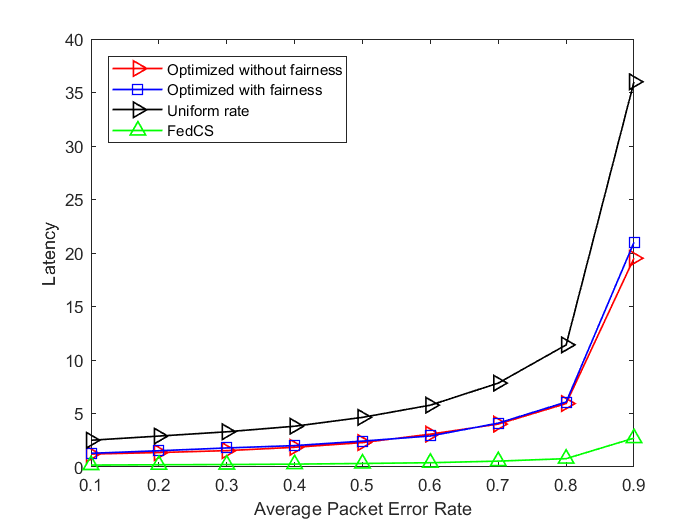}}
\caption{End-to-end delay in the network versus packet error rate for $N=10$, $M=K=5$, and $T=100$. }
\label{fig:delay}
\end{center}
\end{figure}

\begin{figure}[ht]
\begin{center}
\setlength\belowcaptionskip{-2\baselineskip}
\centerline{
\includegraphics[width=0.95\columnwidth]{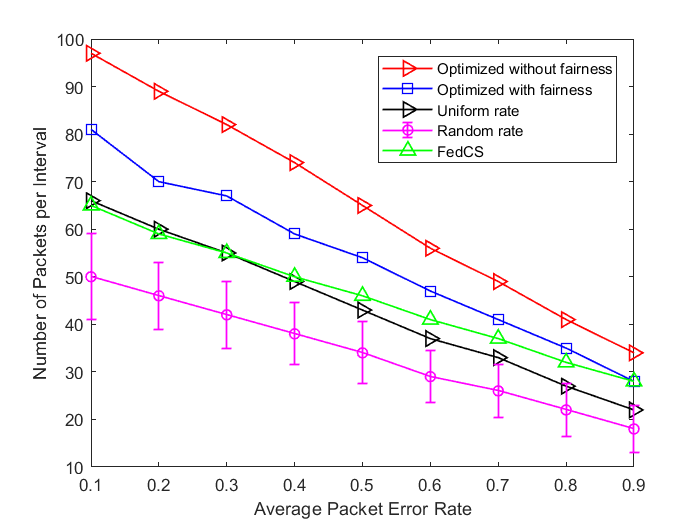}}
\caption{Total throughput of the network (with  $N=10$, $M=K=5$, and $T=100$) using different scheduling policies versus packet error rate. }
\label{fig:th}
\end{center}
\end{figure} 
We investigate scenarios where data samples are split unequally (in terms of class labels) among RSUs to represent \textit{unbalanced} datasets. We evaluated the following methods for each test scenario: i) optimized scheduling without fairness by setting $w_3=0$ in Eqs (\ref{eq:opt}) and (\ref{eq:opt2}) to exclude fairness, ii) optimized scheduling with fairness, where we set $w_1=w_2=w_3$ to weight delay, throughput, and fairness equally (after applying proper normalization for each metric to be in the same scale with zero-mean and unit variance), iii) uniform rate, where all RSUs utilize an equal number of resources, and we set $\alpha_1=\alpha_2=\dots=\alpha_N=M/N$, iv) random rate, where $K=5$ out of $N=10$ RSUs are selected in random to transmit their packets with attempt probability $\alpha_i=M/K$, v) FedCS as an exemplary communication-based method. 
This method requests random clients for their transmission resource information; then, the operator estimates the time required for the distribution and scheduled update and upload steps for a federated learning application. Then, it selects clients that minimize the overall delay \cite{nishio2019client}. 

Fig. \ref{fig:delay} presents the average delay of the system using different methods under different Packet Drop Rates (PDR). As expected, the overall delay increase with PDR, and FedCS yields the best 
performance since it considers delay as its sole performance metric. Nonetheless, the optimized scheduling with and without enforcing diversity outperforms the uniform scheduling.

In Fig. \ref{fig:th}, we can see that the number of successfully transmitted packets per interval drops with PDR due to re-transmission. The decline is linear as (\ref{eq:th1}) suggests. Also, our approach without enforcing fairness ($w_3=0$) results in the highest throughput. This is expected since it regulates scheduling based on the communication performance metrics (delay and throughput); therefore, it outperforms uniform and random scheduling, as well as FedCS (which merely considers delay). It can be seen that enforcing diversity through fairness among class labels declines the system's total throughput. However, we will show that this reduction in throughput will not negatively affect the ultimate learning quality; rather it will increase the learning quality under unbalanced datasets among RSUs.


\begin{figure}
\begin{center}
\setlength\belowcaptionskip{-2\baselineskip}
\centerline{
\includegraphics[width=0.9\columnwidth]{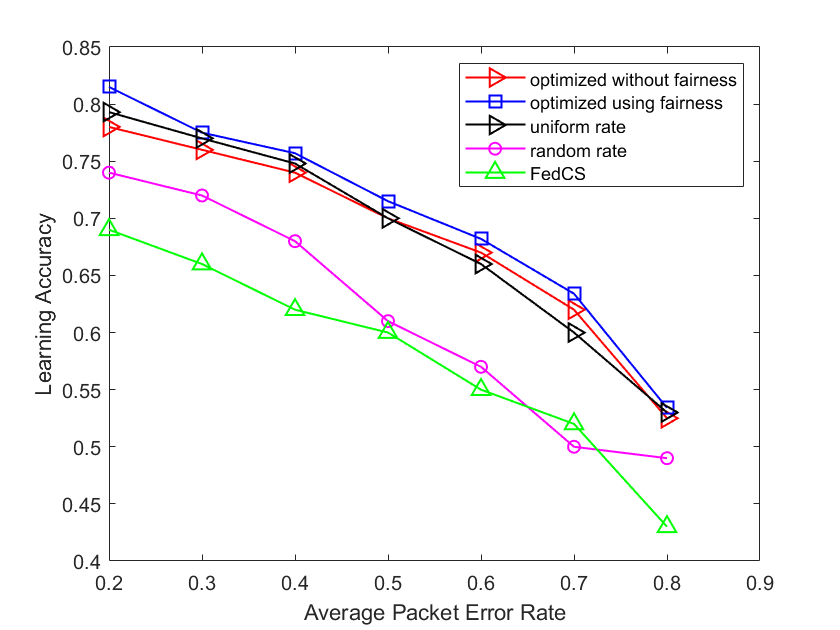}}
\caption{Learning accuracy of the system using different scheduling policies versus packet error rate. The results are for unbalanced CIFAR-10 dataset among RSUs. We set $N=10$, $M=K=5$, and $T=100$. }
\label{fig:cf1}
\end{center}
\end{figure}

\begin{figure}
\begin{center}
\setlength\belowcaptionskip{-2\baselineskip}
\centerline{
\includegraphics[width=0.9\columnwidth]{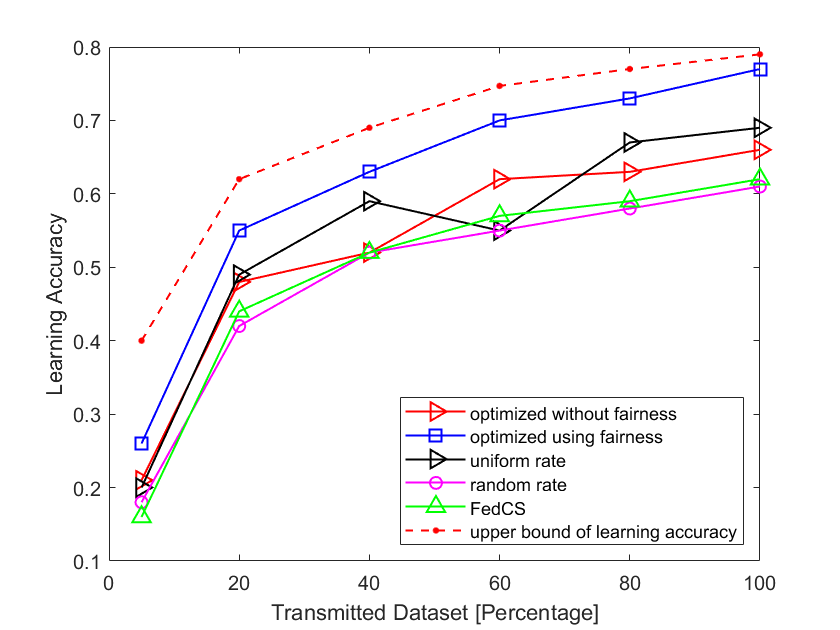}}
\caption{Online learning accuracy of the system using different scheduling policies for different portions of transmitted samples. The results are for an unbalanced CIFAR-10 dataset among RSUs. We set $N=10$, $M=K=5$, and $T=100$. }
\label{fig:cf2}
\end{center}
\end{figure}

It is noteworthy that maximizing throughput is not our primary goal. Indeed, our objective is to maximize the learning quality of the system based on the arrival packets using the proposed scheduling method partially enforced by the uncertainty criterion based on the min-margin metric. To investigate this matter, we evaluate our method by training a CNN on CIFAR10 dataset data samples sent by RSUs under different scheduling policies. 
We consider unbalanced datasets, where ten classes split inequality among RSUs. Indeed, to simulate an unbalanced dataset, we include only 2 out of 10 classes in each RSU dataset. The results are shown in Figs. \ref{fig:cf1} and \ref{fig:cf2}.

It is seen in Fig. \ref{fig:cf1} that the ultimate learning accuracy gradually increases by sending more packets consistently for the optimized scheduling by enforcing fairness on unbalanced datasets. Fig. \ref{fig:cf1} demonstrates a decline in the accuracy by increasing the PDR. In this figure, the proposed optimization when enforcing fairness outperforms all methods for an unbalanced dataset (Figs. \ref{fig:cf1} and \ref{fig:cf2}), since it increases the throughput and the optimized scheduling with enforcing diversity yields superior performance. This highlights the fact that optimizing the scheduling policy merely based on the networking parameters is not optimal for learning-based applications.



\begin{table}[htbp]
\caption{Evaluation with traffic dataset. Three RSUs are used with $\beta_1 = 0.1, \lambda_1 = 1.3$, $\beta_2 = 0.22, \lambda_2 = 1.5$, $\beta_3 = 0.44, \lambda_3 = 1.1$.} 
\begin{center}
\begin{tabular}{ |p{3cm}||p{1.5cm}|p{1.2cm}|p{1.2cm}|  }
 \hline
 & Channel Utilization &Throughput &Learning Quality\\
 \hline
  & RSU1 = 0.58 & &\\
   Optimize without fairness & RSU2 = 0.3 & 89 & 0.85\\
   & RSU3 = 0.12 & &\\
 \hline
  & RSU1 = 0.42 & &\\
   Optimize with fairness & RSU2 = 0.32 & 70 & 0.93\\
   & RSU3 = 0.21 & &\\
 \hline
  & RSU1 = 0.33 & &\\
Uniform rate & RSU2 = 0.33 & 60 & 0.82\\
   & RSU3 = 0.33 & &\\
 \hline 
  & RSU1 = 0.28 & &\\
   Random rate & RSU2 = 0.44 & 46 & 0.81\\
   & RSU3 = 0.17 & &\\
 \hline
   & RSU1 = 0.45 & &\\
   FedCS & RSU2 = 0.5 & 59 & 0.78\\
   & RSU3 = 0.07 & &\\
 \hline
 
 \hline
\end{tabular}
\label{tab1}
\end{center}
\end{table}

Similar results are provided for the traffic light dataset in Table II. Here, each RSU has an unbalanced number of samples in its dataset. Channel utilization rates ($\alpha_1,\alpha_2,\alpha_3$) are given for the RSUs under different policies. It can be observed that FedCS favors RSU2 for its lower expected delay $1/\lambda_2=1/1.5$. On the other hand, the optimized method favors RSU1 for its lower packet drop rate $\beta_1=0.1$. Nevertheless, the attempt probabilities are slightly different with and without enforcing diversity. The optimized scheduling with no diversity maximizes the system throughput whereas the optimized policy with enforcing diversity through fairness among the class labels of the received samples results in the highest learning accuracy of $87\%$.

\section{Conclusion}
In this paper, we investigated the importance of enforcing diversity among collected data samples from RSUs for traffic monitoring applications. We observed that optimal scheduling policies that merely consider networking factors (such as channel drop rate, average delay, etc.), despite maximizing the average throughput and networking efficiency, do not deem optimal in maximizing the ultimate learning accuracy. We examined this condition by simulating unbalanced datasets among RSUs. We offered a new coalition-based greedy optimization that enforces the diversity of the received dataset by imposing fairness on the collected class labels in an interval-by-interval fashion. Then, we use the min-margin criterion to select samples from each class that are less consistent with the trained learning system (hence contributing the most to improving it). Our method outperforms random scheduling, uniform scheduling, and communication-based scheduling methods by a significant margin in terms of learning quality (more than 5\% improvement in classification rate).

\section{Acknowledgment}

This material is based upon the work supported by the National Science Foundation under Grant Numbers 2008784 and 2204721.

\bibliographystyle{ieeetr}
\bibliography{main} 

\end{document}